\documentstyle[11pt,epsfig,epsf]{article}
\input{amssymb.sty}

\begin{document}

\title{\bf Fredholm Perturbation of Spectra of $2\times 2$ Upper Triangular Matrix  \footnote{This work is
supported by the NSF of China (Grant Nos. 10771034, 10771191 and
10471124) and the NSF of Fujian Province of China (Grant Nos.
Z0511019, S0650009).}}
\author{{Shifang Zhang$^{1,2}$, \, Huaijie Zhong$^2$, \, Junde Wu$^1$\footnote{Corresponding
author: Junde Wu, E-mail: wjd@zju.edu.cn}}
\\\\$^1$\small\it Department of Mathematics, Zhejiang University,
Hangzhou 310027, P. R. China
\\\\$^2$\small\it Department of Mathematics, Fujian Normal University, Fuzhou 350007, P. R. China, \\shifangzhangfj@163.com, hjzhong@fjnu.edu.cn}
 \date{}\maketitle
\begin{center}
\begin{minipage}{140mm}
{\bf Abstract} { \, 
As we knew, study the perturbation theory of spectra of operator is a very important project in mathematics physics, in particular, in quantum mechanics. In this paper, we characterize the Fredholm perturbation for the Weyl spectrum, essential spectrum, spectrum,
left spectrum, right spectrum, lower semi-Fredholm spectrum, upper semi-Weyl spectrum and lower semi-Weyl spectrum of upper triangular operator matrix $M_{C}=\left( \begin{array}{cc} A&C\\ 0&B\\ \end{array}\right)$.}

\vspace{4mm} {\bf Keywords}\,\,\,{Operator matrix;
spectra; perturbation.}

 \vspace{4mm} {\bf  AMS classifications}\,\,\,{47A10} 

\end{minipage}
\end{center}
\vskip 0.2in
%
%
 \section{\Large\bf ~~Introduction}
 \vskip 0.1in

Let $H$ and $K$ be the complex infinite
dimensional separable Hilbert spaces, $B(H, K)$ be the set of
all bounded linear operators from $H$ into $K$. For simplicity, we
write $B(H, H)$ as $B(H).$ If $T\in B(H, K)$, we use $R(T)$ and $N(T)$ to denote
the range and kernel of $T$, respectively, and define $\alpha(T)=\dim N(T)$
and $\beta (T)=\dim (K/ R(T))$. For $T\in B(H, K)$, if $R(T)$ is
closed and $\alpha (T)<\infty$, we call $T$ an upper semi-Fredholm
operator; if $\beta (T)<\infty$, then $T$ is called a lower
semi-Fredholm operator. If $T$ is either an upper or lower
semi-Fredholm operator, then $T$ is called a semi-Fredholm operator. In this case, the index of $T$ is defined as
ind$(T)=\alpha(T)-\beta(T).$ 
If $T$ is a semi-Fredholm operator with $\alpha (T)<\infty$ and $\beta (T)<\infty$, then $T$ is called a Fredholm operator.
 For $T\in B(H)$, the ascent asc$(T)$ and the descent des$(T)$ are given by 
asc$(T)=\inf\{k \geqslant 0: N(T^{k})=N(T^{k+1})\}$ and des$(T)=\inf\{k \geqslant 0:R(T^{k})=R(T^{k+1})\}$, respectively;
the infimum over the empty set is taken to be $\infty$.

\vspace{3mm}

Let $G(H, K), G_l(H, K)$, $G_r(H, K), \Phi(H, K), \Phi_+(H, K)$ and $\Phi_-(H, K)$, respectively,
 denote the sets of all invertible operators, left invertible operators, right
 invertible operators, Fredholm operators, upper semi-Fredholm operators and lower
semi-Fredholm operators from $H$ into $K$. The sets of all Weyl operators, upper semi-Weyl operators and lower
semi-Weyl operators from $H$ into $K$ are defined, respectively, by

\vspace{3mm}
\par  $\Phi_{0}(H, K):=\{T \in \Phi(H, K):$ ind$(T)=0 \},$
\par  $\Phi_+^-(H, K):=\{T \in \Phi_+(H, K):$ ind$(T)\leq 0\},$
\par  $\Phi_-^+(H, K):=\{T \in \Phi_-(H, K):$ ind$(T)\geq 0\}.$

 \noindent When $H=K$, the above 9 kind operator classes are also abbreviated as  $G(H), G_l(H)$, $G_r(H),
\Phi(H),$ $ \Phi_+(H)$, $\Phi_-(H),$ $\Phi_0(H)$,$\Phi_+^-(H)$ and $\Phi_-^+(H),$ respectively.

For $T\in B(H)$, its corresponding spectra are, respectively, defined by

\vspace{3mm}
   \par the spectrum: $\sigma_{}(T)= \{\lambda \in {\mathbb{C}}:T-\lambda I \makebox{ is not invertible}\},$
  \par the left spectrum: $\sigma_{l}(T)=\{\lambda \in {\mathbb{C}}:T-\lambda I$ is not left invertible$\},$
 \par the right spectrum: $\sigma_{r}(T)=\{\lambda \in {\mathbb{C}}:T-\lambda I$ is not right invertible$\},$
 \par the essential spectrum:  $\sigma_{e}(T)=\{\lambda \in {\mathbb{C}}: T-\lambda I\not \in \Phi(H)\},$
  \par the upper semi-Fredholm spectrum: $\sigma_{SF+}(T)=\{\lambda \in {\mathbb{C}}: T-\lambda I \not \in \Phi_+(H)\},$
  \par the lower semi-Fredholm spectrum:  $\sigma_{SF-}(T)=\{\lambda \in {\mathbb{C}}: T-\lambda I \not \in \Phi_-(H)\},$
  \par the Weyl spectrum: $\sigma_{w}(T)=\{\lambda \in {\mathbb{C}}:T-\lambda I \not \in \Phi_{0}(H) \},$
  \par the upper semi-Weyl spectrum: $\sigma_{aw}(T)=\{\lambda \in {\mathbb{C}}: T-\lambda I \not \in \Phi_+^-(H)\},$
  \par the lower semi-Weyl spectrum: $\sigma_{sw}(T)=\{\lambda \in {\mathbb{C}}: T-\lambda I \not \in \Phi_{-}^+(H)\},$
  \par the Browder spectrum: $\sigma_{b}(T)=\{\lambda \in {\mathbb{C}}: T-\lambda I \not \in \Phi_{b}(H)\},$
  where  $\Phi_{b}(H):=\{T \in \Phi(H):$ asc$(T)<\infty$ and  des$(T)<\infty\}.$
 
\vspace{3mm}

It is well known that all the above spectra are compact
nonempty subsets of complex plane ${\mathbb{C}}$.

\vspace{3mm}

Let $H$ be a Hilbert space and $T$ be a bounded linear operator defined on $H$ and $H_1$ be 
an invariant closed subspace of $T$. Then $T$ can be represented by 
the form of $$T=\left(
       \begin{array}{cc}
        *&*\\
         0&*\\
       \end{array}
     \right):H_1\oplus H_1^{\perp}\rightarrow H_1\oplus H_1^{\perp},$$
which motivated the interest in $2 \times 2$ upper-triangular
operator matrices (see [1-19]). 

Henceforth, for $A\in B(H)$, $B\in B(K)$ and $C\in B(K, H)$, we put
$M_{C}=\left( \begin{array}{cc} A&C\\ 0&B\\ \end{array}\right)$. It is clear that $M_C\in B(H\oplus K)$. 
Recent, people studied the perturbation theory of some spectra of $M_{C}$, for example, in [8], for the spectrum $\sigma (M_C)$,
the perturbation result is $$\bigcap_{C\in B(K,\,H)}\sigma(M_{C})=\sigma_{l} (A)\cup \sigma_{r}(B) \cup\{\lambda\in {\mathbb{C}}:
\alpha(B-\lambda)\not=\beta(A-\lambda)\}.\eqno (1)$$ In [5], for the Weyl spectrum $\sigma_w (M_C)$ and the essential spectrum
$\sigma_e (M_C)$, the perturbation results are $$\bigcap_{C\in B(K,\,H)}\sigma_{w}(M_{C})=\sigma_{SF+} (A)\cup \sigma_{SF-}(B)
\cup\{\lambda\in {\mathbb{C}}:\alpha(A-\lambda)+\alpha(B-\lambda)\not=\beta(A-\lambda)+\beta(B-\lambda)\}\eqno (2)$$ and
$$\bigcap_{C\in B(K,\,H)}\sigma_{e}(M_{C})=\sigma_{SF+} (A)\cup \sigma_{SF-}(B)\cup$$$$\{\lambda\in {\mathbb{C}}:
\min(\beta(A-\lambda),\alpha(B-\lambda))<\max(\beta(A-\lambda),\alpha(B-\lambda))=\infty\}.\eqno (3)$$
 In [1-3, 10], the authors also characterize completely sets $\bigcap_{C\in B(K,\,H)}\sigma_{*}(M_{C})$, where 
$\sigma_{*}(M_{C})$ may be the Browder spectrum, left spectrum, right spectrum, lower semi-Fredholm spectrum, upper semi-Fredholm spectrum,  
lower semi-Weyl spectrum or upper semi-Weyl spectrum of $M_{C}$, respectively.

Moreover, in [13-15], for the spectra $\sigma_{*}(M_{C})$, where $\sigma_{*}=\sigma_{r},\sigma_{SF-}$ or $\sigma_{sw}$, its  perturbation result is $$\bigcap_{C\in
G(K,\,H)}\sigma_{*}(M_{C})=(\bigcap_{C\in B(K,\,H)}\sigma_{*}(M_{C}))\cup \{\lambda\in {\mathbb{C}}:
A-\lambda\, \makebox {is  compact}\}; \eqno (4)$$ for the spectra $\sigma_{*}(M_{C})$, where $\sigma_{*}=\sigma_{l},\sigma_{SF+}$
or $\sigma_{aw}$, its perturbation result is $$\bigcap_{C\in G(K,\,H)}\sigma_{*}(M_{C})=(\bigcap_{C\in
B(K,\,H)}\sigma_{*}(M_{C}))\cup \{\lambda\in {\mathbb{C}}:B-\lambda\, \makebox {is  compact}\};\eqno (5)$$
for the spectra $\sigma_{*}(M_{C})$, where $\sigma_{*}=\sigma,\sigma_{e}$ 
or $\sigma_{w}$, its  perturbation result is $$\bigcap_{C\in G(K,\,H)}\sigma_{*}(M_{C})=(\bigcap_{C\in
B(K,\,H)}\sigma_{*}(M_{C}))\cup \{\lambda\in {\mathbb{C}}: A-\lambda\, \makebox{or}\, B-\lambda\, \makebox {is compact}\}.\eqno (6)$$

Note that equations (1) to (3) showed the perturbation of all bounded linear operator $C$ in $B(K, H)$, 
and equations (4) to (6) showed the perturbation of all bounded invertible linear operator $C$ in $G(K, H)$.

\vspace{3mm}

In this paper, we characterize the Fredholm perturbation for the Weyl spectrum, essential spectrum, spectrum,
left spectrum, right spectrum, lower semi-Fredholm spectrum, upper semi-Weyl spectrum and lower semi-Weyl spectrum of $M_{C}$.

%
%
%
%

\vskip 0.2in \section{\Large\bf ~~Main results and proofs}
 \vskip 0.2in
 
At first, in order to characterize the perturbation of Weyl spectrum of $M_C$, we need the following: 

\noindent{\bf Lemma 1.} For a given pair $(A, B)\in B(H)\times B(K)$, the following statements are equivalent:

(i). there exists some $C\in B(K, H)$ such that $M_C \in \Phi_0(H\oplus K),$

(ii). $A \in \Phi_+(H)$, $B\in \Phi_-(K)$ and $\alpha(A)+\alpha(B)=\beta(A) + \beta(B),$

(iii). there exists some $Q\in G(K, H)$ such that  $M_Q \in\Phi_0(H\oplus K),$

(iv). there exists some $Q\in \Phi(K, H)$ such that  $M_Q \in\Phi_0(H\oplus K).$

\noindent{\bf Proof.} (i)$\Leftrightarrow$(ii) was proved in [5, Theorem 3.6]. 

\vspace{2mm}

(ii)$\Rightarrow$(iii). It is sufficient to prove that if $A \in
\Phi_+ (H)$, $B\in \Phi_- (K)$ and $\beta(A)=\alpha(B)=\infty,$ then
there exists $Q\in G(K,H)$ such that $M_Q \in \Phi_0(H\oplus K).$ To
show this, there are three cases to consider:

\vspace{2mm}

 Case 1. Suppose $\alpha(A)=\beta(B)<\infty$. Define an operator
$Q : K\rightarrow H\, \makebox{ by}\,
Q = \left(
       \begin{array}{cc}
        T_1&0\\
         0&T_2\\
       \end{array}
     \right):
       N(B)\oplus
         N(B)^\perp
      \rightarrow
     R(A)^\perp \oplus R(A),$ where $T_1 $ and $T_2$ are invertible operators. Obviously, $Q\in
G(K,H)$ and $M_Q \in \Phi(H\oplus K)$. Also, it is
evident that $N(M_Q)=N(A)\oplus \{0\}$ and $R(M_Q)^\perp=\{0\}\oplus
R(B)^\perp $. Thus $\alpha(M_Q)=\beta(M_Q)= \alpha(A)=\beta(B)<\infty$, and
hence $ M_Q \in \Phi_0(H\oplus K)$ is clear.

\vspace{2mm}

Case 2. Suppose $\beta(B)<\alpha(A)<\infty$ and put
$l=\alpha(A)-\beta(B).$  Note that $\beta(A) =\dim
N(B)^\perp=\infty$, let $ R(A)^\perp=H_1\oplus H_2$ and $\dim
H_2 =l$, $ N(B)^\perp =K_1\oplus K_2$ and $ \dim(K_1)=l.$

Define an operator
$Q : K\rightarrow H\, \makebox{ by}\,
Q = \left(
       \begin{array}{ccc}
        T_1&0 &0\\
         0&T_2&0\\
   0&0&T_3\\  \end{array}
     \right):
       N(B)\oplus
         K_1\oplus K_2
      \rightarrow
      H_1\oplus H_2 \oplus R(A),$ where $T_1, T_2$ and $T_3$ are invertible operators. Obviously,
$Q\in B(K,H)$ is invertible. Now we claim that $M_Q\in \Phi_0(H\oplus
K)$. In fact, $M_Q$ has the following form:
$M_Q=\left(\begin{array}{ccccc}
   0&  0&T_1&0&0\\0&0&0&T_2&0\\0&A_1&0& 0&T_{3}\\0&0&0&B_{1}&B_{2} \\0&0&0&0&0 \\
 \end{array} \right):{N(A)}\oplus{N(A)}^\perp \oplus{N(B)}\oplus{K_1} \oplus{K_2}  \longrightarrow
    H_1 \oplus{H_2}\oplus{R(A)}\oplus{R(B)}\oplus{R(B)}^\perp,$
where $A_1\in B(N(B)^\perp, R(A))$ and $(B_1\,\, B_2)\in B((K_1\oplus K_2), R(B))$ are invertible operators.
Moreover, observe that $\dim K_1<\infty,$ we have $B_1\in G(K_1, R(B_1))$,
$B_2\in G(K_2, R(B_2))$ and $\dim K_1=\dim R(B_1)=\dim (R(B)\ominus R(B_2)).$

Now let $W_1=\left(\begin{array}{ccccc}
I& 0&0&0&0\\0&I&0&0&0\\0&0&I& 0&0\\0&-B_1T_2^{-1}&0&I&0 \\0&0&0&0&I \\ \end{array}
\right):{N(A)}\oplus{N(A)}^\perp \oplus{N(B)}\oplus{K_1} \oplus{K_2}  \longrightarrow
{N(A)}\oplus{N(A)}^\perp \oplus{N(B)}\oplus{K_1} \oplus{K_2},$

and $W_2=\left(    \begin{array}{ccccc}
I& 0&0&0&0\\0&I&0&0&-A_1^{-1}T_3\\0&0&I& 0&0\\0&0&0&I&0 \\0&0&0&0&I \\ \end{array}
\right):H_1 \oplus{H_2}\oplus{R(A)}\oplus{R(B)}\oplus{R(B)}^\perp \longrightarrow
H_1 \oplus{H_2}\oplus{R(A)}\oplus{R(B)}\oplus{R(B)}^\perp.$

\noindent Then $W_1M_QW_2=\left( \begin{array}{ccccc}
0& 0&T_1&0&0\\0&0&0&T_2&0\\0&A_1&0& 0&0\\0&0&0&0&B_{2} \\0&0&0&0&0 \\  \end{array}
\right):{N(A)}\oplus{N(A)}^\perp \oplus{N(B)}\oplus{K_1} \oplus{K_2}  \longrightarrow
 H_1 \oplus{H_2}\oplus{R(A)}\oplus{R(B)}\oplus{R(B)}^\perp.$
Since $A_1, T_1 $ and $T_2 $ are invertible, we get that 
 $R(W_1M_QW_2)= H_1\oplus{H_2}\oplus{R(A)}\oplus{R(B_2)}\oplus\{0\}$ and
$N(W_1M_QW_2)= N(A)\oplus\{0\}\oplus\{0\}\oplus\{0\}\oplus\{0\},$ and 
$R(W_1M_QW_2)^\perp=\{0\}\oplus\{0\}\oplus\{0\}\oplus({R(B)}\ominus{R(B_2)})\oplus{R(B)^\perp}.$
Thus $W_1M_QW_2 \in\Phi(H\oplus K)$ and 

\begin{eqnarray*}
& \alpha (W_1M_QW_2) &=\alpha(A)=l+\beta(B)\\
&&=\dim K_1+\beta(B)\\
& &=\dim ({R(B)}\ominus{R(B_2)})+\beta(B)\\
& &=\beta(W_1M_QW_2)<\infty.\\
\end{eqnarray*}

So $W_1M_QW_2 \in\Phi_0(H\oplus K).$ Also since $W_1$ and $W_2$ are
invertible, it follows that $M_Q \in\Phi_0(H\oplus K).$

\vspace{2mm}

Case 3. Suppose $\alpha(A)<\beta(B)<\infty$, put
$l=\beta(B)-\alpha(A).$ Since $\dim R(A)=\dim N(B)=\infty$, let
R(A)$=H_1\oplus H_2$ and $\dim H_1 =l$, $ N(B)=K_1\oplus K_2$
and $\dim(K_2)=l.$ That $\dim H_2=\dim(K_1)=\infty$ is clear. Define an operator
$Q : K\rightarrow H\, \makebox{ by}\,
Q = \left(
       \begin{array}{ccc}
        T_1&0 &0\\
         0&T_2&0\\
   0&0&T_3\\
       \end{array}
     \right):
         K_1\oplus K_2 \oplus N(B)^\perp
      \rightarrow
      R(A)^\perp\oplus  H_1\oplus H_2,$ where $T_1,$ $T_2$ and $T_3$ are invertible operators. Obviously,
$Q\in G(K,H)$. Similar to the proof of Case 2, we can
also show that $M_Q \in \Phi_0(H\oplus K)$. 

It follows from Case 1 to Case 3 that (ii)$\Rightarrow$(iii). 

 Finally, (iii)$\Rightarrow$(iv) and (iv)$\Rightarrow$(i) are clear. The lemma is proved.
 
 \vspace{3mm}

From Lemma 1 and Equation (2), we have the following:

\vspace{3mm}

\noindent{\bf Theorem 1.} For a given pair $(A, B)\in B(H)\times B(K)$, we have
$$\bigcap_{C\in
\Phi(K,\,H)}\sigma_{w}(M_{C})=\bigcap_{C\in G(K,\,H)}\sigma_{w}(M_{C})=\bigcap_{C\in
B(K,\,H)}\sigma_{w}(M_{C})$$$$=\sigma_{SF+} (A)\cup \sigma_{SF-}(B)\cup\{\lambda\in {\mathbb{C}}:
\alpha(A-\lambda)+\alpha(B-\lambda)\not=\beta(A-\lambda)+\beta(B-\lambda)\}.$$

In order to characterize the perturbation of essential spectrum of $M_C$, we need the following:

\noindent{\bf Lemma 2.} For a given pair $(A, B)\in B(H)\times B(K)$, the following statements are equivalent:
 \par(i). there exists some  $C\in B(K, H)$ such that $M_C \in \Phi(H\oplus K),$
 \par(ii). $\left\{\begin{array}{l}  A \in \Phi (H)\,\,\makebox{and}\,\, B \in \Phi (K)
      \\\makebox{or}\,\, A\in \Phi_+(H), B\in \Phi_-(K)\,\, \makebox{and}\,\,
  \beta(A)=\alpha(B)=\infty,\end{array}\right.$
 \par(iii).  there exists some $Q\in G(K, H )$ such that $M_Q \in \Phi(H\oplus K),$
\par(iv). there exists some $Q\in \Phi(K, H )$ such that $M_Q \in \Phi(H\oplus K).$

\vspace{3mm}

\noindent{\bf Proof.} (i)$\Rightarrow$(ii).  Suppose that $M_C \in \Phi(H\oplus K)$ for some $C\in B(K, H)$. It follows from 
[5, Theorem 3.2] that $A\in \Phi_+(H), B\in \Phi_-(K)$. Moreover, by [19, Lemma 2.2] we have that either both $A$ and $B$ are Fredholm operators 
or neither $A$ nor $B$ is a Fredholm operator. Thus $\beta(A)=\alpha(B)=\infty$ when neither $A$ nor $B$ is a  Fredholm operator.

(ii)$\Rightarrow$(iii). To do this, if $A\in \Phi(H)$ and $B\in \Phi(K)$, then $M_C\in \Phi (H\oplus K)$ 
for every $C\in B(K, H )$. On the other hand, if  $A\in \Phi_+(H)$, $B\in \Phi_- (K)$ and
$\beta(A)=\alpha(B)=\infty.$ Define an operator
$Q : K\rightarrow H\, \makebox{ by}\,
Q=\left(\begin{array}{cc} T_1&0\\ 0&T_2\\ \end{array} \right):
N(B)\oplus N(B)^\perp \rightarrow R(A)^\perp \oplus R(A)$,
where $T_1$ and $T_2$ are invertible operators. Obviously,
$Q\in G(K, H)$, and it is easy to show that $M_Q\in\Phi(H\oplus K)$.

(iii) $\Rightarrow$ (iv) and (iv) $\Rightarrow$ (i)
are obvious. The lemma is proved.

\vspace{3mm}

 From Lemma 2 and Equation (3) we have the following immediately:

\noindent{\bf Theorem 2.} For a given pair $(A, B)\in B(H)\times B(K)$, we have
$$\bigcap_{C\in \Phi(K,\,H)}\sigma_{e}(M_{C})=\bigcap_{C\in G(K,\,H)}\sigma_{e}(M_{C})
=\bigcap_{C\in B(K,\,H)}\sigma_{e}(M_{C})$$$$=\sigma_{SF+} (A)\cup \sigma_{SF-}(B)\cup\{\lambda\in {\mathbb{C}}:
\min(\beta(A-\lambda),\alpha(B-\lambda))<\max(\beta(A-\lambda),\alpha(B-\lambda))=\infty\}.$$

\vspace{3mm}

In order to characterize the perturbation of spectrum of $M_C$, we need the following lemma 
which is a generalization in [9, Theorem 2] in the case of Hilbert spaces:

\vspace{3mm}

\noindent{\bf Lemma 3.} For a given pair $(A, B)\in B(H)\times B(K)$, the following statements are equivalent:
 \par(i). there exists some $C\in B(K, H )$ such that $M_C$ is invertible,
 \par(ii). $A$ is left invertible, $B$ is right invertible and $\beta(A)=\alpha(B),$
 \par(iii). there exists some $Q\in G(K, H )$ such that $M_Q$ is invertible,
 \par(iv). there exists some $Q\in \Phi(K, H )$ such that $M_Q$ is invertible.

\noindent{\bf Proof.} (i)$\Rightarrow$(ii) is prove in [9, Theorem 2]. In fact, if $M_C$ is invertible, it is easy to show that 
$A$ is left invertible and $B$ is right invertible, which implies that $\alpha(A)=\beta(B)=0$. Moreover, it follows from Lemma 1 that
$\alpha(A)+\alpha(B)=\beta(A)+\beta(B),$ thus $\beta(A)=\alpha(B)$.

(ii)$\Rightarrow$(iii). Suppose $A$ is left invertible, $B$ is right invertible and $\beta(A)=\alpha(B)$. Define an operator $Q : K\rightarrow H\, \makebox{ by}\,
Q=\left(\begin{array}{cc} T_1&0\\
 0&T_2\\\end{array}\right):N(B)\oplus
         N(B)^\perp\rightarrow
R(A)^\perp \oplus R(A),$ where $T_1 $ and $T_2$ are invertible operators. it is evident that  $Q\in G(K,H)$ and $M_Q \in G(H\oplus K)$.

(iii) $\Rightarrow$ (iv) and (iv) $\Rightarrow$ (i)
are obvious. The lemma is proved. 

\vspace{3mm}

From Lemma 3 and Equation (1), the following theorem is immediate:

\vspace{3mm}

\noindent{\bf Theorem 3.} For a given pair $(A, B)\in B(H)\times B(K)$, We have
 $$\bigcap_{C\in \Phi(K,\,H)}\sigma(M_{C})=\bigcap_{C\in G(K,\,H)}\sigma(M_{C})
 =\bigcap_{C\in B(K,\,H)}\sigma(M_{C})$$$$=\sigma_{l} (A)\cup \sigma_{r}(B) \cup\{\lambda\in {\mathbb{C}}:
\alpha(B-\lambda)\not=\beta(A-\lambda)\}.$$

\vspace{3mm}

In order to characterize the perturbation for left spectrum, right spectrum, lower semi-Weyl spectrum, upper semi-Weyl spectrum  
and lower semi-Fredholm spectrum of $M_C$, we need the following three lemmas:
\vspace{3mm}

\noindent{\bf Lemma 4.} For a given pair $(A, B)\in B(H)\times B(K)$, if either $A$ or $B$ is a
compact operator, then for each $C\in\Phi(K ,H)$, $M_C$ is not a semi-Fredholm operator.

\noindent{\bf Proof.} If $B$ is a compact operator,  then we can claim that $M_C$ is not a semi-Fredholm
operator for each $C\in \Phi(K, H)$.  If not, assume that $C_0\in \Phi(K ,H)$ such that
$M_{C_0}$ is a semi-Fredholm operator. Since $C_0\in \Phi(K ,H)$,
there exists $C_1\in \Phi(H ,K)$ such that $C_0C_1=I+K$, where $K\in B(H)$ is a compact
operator. Note that $$\left(\begin{array}{cc}A&C_0\\0&B\\
\end{array} \right)\left(\begin{array}{cc}
I&0\\ -C_1A&I\\ \end{array}\right)=\left(\begin{array}{cc}
A-C_0C_1A&C_0\\ -BC_1A&B\\ \end{array}\right)
=\left(\begin{array}{cc} -KA&C_0\\ -BC_1A&B\\\end{array}\right),$$
we have that $\left(\begin{array}{cc}-KA&C_0\\
-BC_1A&B\\\end{array}\right)$ is a semi-Fredholm operator. Also since $K$
and $B$ are compact operators,   both $ \left(\begin{array}{cc}
 0&0\\-BC_1A&0\\\end{array}\right)$ and $\left(\begin{array}{cc}
-KA&0\\0&B\\\end{array}\right)$ are also compact. Thus
$\left(\begin{array}{cc}0&C_0\\ 0&0\\\end{array}\right)$ is a
semi-Fredholm operator, which is impossible. 
So $M_C$ is not a semi-Fredholm operator for each $C\in \Phi(K, H)$.

Similarly, we can prove when $A$ is a compact operator, $M_C$ is not a semi-Fredholm
operator for each $C\in \Phi(K, H)$. The lemma is proved.

\vspace{3mm}

\noindent{\bf Lemma 5.} The following statements are equivalent:
\par(i). $B$ is not compact,
\par(ii). for each given $A\in\Phi_+(H)$, if $\beta(A) =\infty$, then there exists an operator $C\in G(K, H)$  such that $M_C $ is an upper semi-Weyl operator and $\alpha(M_C)=\alpha(A)$,
\par(iii). for each given $A\in\Phi_+(H)$, if $\beta(A) =\infty$, then there exists an operator $C\in \Phi(K, H)$ such that $M_C $ is an upper semi-Weyl operator and $\alpha(M_C)=\alpha(A)$,
\par(iv). for each given $A\in\Phi_+(H)$, if $\beta(A) =\infty$, then there exists an operator $C\in G(K, H)$ such that $M_C$ is an upper semi-Weyl operator,
\par(v). for each given $A\in\Phi_+(H)$, if $\beta(A) =\infty$, then there exists an operator $C\in\Phi(K, H )$ such that $M_C$ is an upper semi-Weyl operator,
\par(vi). for each given $A\in\Phi_+(H)$, if h $\beta(A) =\infty$, then there exists an operator $C\in G(K, H)$ such that $M_C$ is an upper semi-Fredholm operator,
\par(vii). for each given $A\in\Phi_+(H)$, if $\beta(A) =\infty$, then there exists an operator $C\in \Phi(K, H)$ such that $M_C$ is an upper semi-Fredholm operator.

\noindent{\bf Proof.}  Obviously, we only need to prove the implications (i) $\Rightarrow$ (ii) and (vii) $\Rightarrow$ (i).

 (vii) $\Rightarrow$ (i). If $B$ is compact, then it follows from
Lemma 4 that $M_C$ is not a semi-Fredholm operator for each $C\in
\Phi(K, H)$, which contradicts with (vii). Thus $B$ is not compact.

 (i) $\Rightarrow$ (ii). Suppose that $B$ is not compact.
Then we consider the following two cases:

Case 1. Assume that $R(B)$ is closed. Since the assumption that $B$
is not compact, we have that $\dim{N(B)}^\perp=\infty.$ Also since
$\beta(A)=\infty$, let ${R(A)}^\perp= H_1 \oplus H_2$ with $\dim H_1
=\dim N(B)$ and $\dim H_2 =\infty.$ Define an operator $C :
K\rightarrow H$ by $$C=\left(\begin{array}{cc}C_1&0\\0&C_2\\\end{array}
 \right):{N(B)}\oplus{N(B)}^\perp\longrightarrow {H_1}\oplus(H_2\oplus R(A)),$$
where $C_1\in B(N(B),H_1)$ and $C_2\in B({N(B)}^\perp,H_2\oplus
R(A))$ are invertible operators. Obviously, operator $C$ is invertible. By [12, Lemma 2],
$M_C$ is an upper semi-Fredholm operator. Moreover, it is easy to prove that 
$N( M_C)=N(A)\oplus\{0\}$ and $\dim {R( M_C)}^\perp\geq\dim H_2=\infty.$ Thus, $M_C$ is an upper
semi-Weyl operator and $\alpha (M_C)=\alpha (A)$.

Case 2. Assume that $R(B)$ is not closed. By [13, Lemma 3.6] and its proof, we can obtain an operator
 $C\in G(K, H)$ such that $M_C$ is an upper semi-Weyl operator and $\alpha (M_C)=\alpha (A)$. The lemma is proved. 

\vspace{3mm}

Duality, we have:

\vspace{3mm}
 
\noindent{\bf Lemma 6.} The following statements are equivalent:
\par(i). $A$ is not compact,
\par(ii). for each given $B\in\Phi_-(K)$, if $\alpha(B)=\infty$, then there exists an operator $C\in G(K, H)$  such that $M_C $ is a lower semi-Weyl operator and $\beta(M_C)=\beta(B)$,
\par(iii). for each given $B\in\Phi_-(K)$, if $\alpha(B)=\infty$, then there exists an operator $C\in \Phi(K, H)$ such that $M_C $ is a lower semi-Weyl operator and $\beta(M_C)=\beta(B)$,
\par(iv). for each given $B\in\Phi_-(K)$, if $\alpha(B)=\infty$, then there exists an operator $C\in G(K, H)$ such that $M_C$ is a lower semi-Weyl operator,
\par(v). for each given  $B\in\Phi_-(K)$, if $\alpha(B)=\infty$, then there exists an operator $C\in\Phi(K, H )$ such that $M_C$ is a lower semi-Weyl operator,
\par(vi). for each given $B\in\Phi_-(K)$, if $\alpha(B)=\infty$, then there exists an operator $C\in G(K, H)$ such that $M_C$ is a lower semi-Fredholm operator,
\par(vii). for each given $B\in\Phi_-(K)$, if $\alpha(B)=\infty$, then there exists an operator $C\in \Phi(K, H)$ such that $M_C$ is a lower semi-Fredholm operator.

\vspace{3mm}

Our Theorem 4 and Theorem 5 following show the similar conclusions as Equation (4)-(5). 

\vspace{3mm}

\noindent{\bf Theorem 4.} For a given pair $(A, B)\in B(H)\times B(K)$, we have $$\bigcap_{C\in
\Phi(K,\,H)}\sigma_{*}(M_{C})=(\bigcap_{C\in
B(K,\,H)}\sigma_{*}(M_{C}))\cup \{\lambda\in {\mathbb{C}}:
A-\lambda\, \makebox{is  compact}\},$$ where
$\sigma_{*} \in\{\sigma_{r},\sigma_{SF-},\sigma_{sw}\}.$

\noindent{\bf Proof.}  According to Lemma 4, it is clear that
$$\bigcap_{C\in \Phi(K,\,H)}
\sigma_{*}(M_{C})\supseteq (\bigcap_{C\in B(K,\,H)}\sigma_{*}(M_{C}))\cup \{\lambda\in
 {\mathbb{C}}: A-\lambda\, \makebox{ is compact}\}.$$

In order to show the theorem, we only need to prove that
$$\bigcap_{C\in\Phi(K,\,H)}\sigma_{*}(M_{C})\subseteq (\bigcap_{C\in
B(K,\,H)}\sigma_{*}(M_{C}))\cup \{\lambda\in {\mathbb{C}}:
A-\lambda\, \makebox{ is compact}\}.$$

(i). Suppose that $\sigma_{*}(\cdot)=\sigma_{SF-}(\cdot)$ and 
$\lambda \not\in (\bigcap_{C\in B(K,\,H)}\sigma_{SF-}(M_{C}))\cup
\{\lambda\in {\mathbb{C}}: A-\lambda\, \makebox{ is compact}\}.$
 Then $A-\lambda$ is not compact and there exists $C\in B(K,\,H)$ such
 that $M_{C}-\lambda\in \Phi_-(H\oplus K),$  and hence  $B-\lambda\in \Phi_-(K).$

Case 1. $ \alpha(B-\lambda) =\infty$. It follows from Lemma 6 that there exists 
$C\in\Phi(K, H)$ such that $M_C-\lambda$ is a lower semi-Fredholm operator.  
This implies that $\lambda \not\in \bigcap_{C\in \Phi(K,\,H)}\sigma_{SF-}(M_{C})$.
It is clear that $$\bigcap_{C\in \Phi(K,\,H)}\sigma_{SF-}(M_{C})
\subseteq (\bigcap_{C\in B(K,\,H)}\sigma_{SF-}(M_{C}))\cup \{\lambda\in
{\mathbb{C}}: A-\lambda\, \makebox{ is compact}\}.$$

Case 2. $\alpha(B-\lambda)<\infty$. This implies that $B-
\lambda\in\Phi(K),$ and so $A-\lambda\in\Phi_-(H)$ since $M_{C}-\lambda\in \Phi_-(H\oplus K).$ Thus, we have
that $M_C-\lambda$ is a lower semi-Fredholm operator for each $C\in
B(K,\,H)$, which means $\lambda \not\in \bigcap_{C\in
\Phi(K,\,H)}\sigma_{SF-}(M_{C})$. Thus
$$\bigcap_{C\in\Phi(K,\,H)}\sigma_{SF-}(M_{C})
\subseteq \bigcap_{C\in B(K,\,H)}\sigma_{SF-}(M_{C})\cup \{\lambda\in
{\mathbb{C}}: A-\lambda\, \makebox{ is compact}\}.$$

Together Case 1 with Case 2, we have 

$$\bigcap_{C\in\Phi(K,\,H)}\sigma_{SF-}(M_{C})=(\bigcap_{C\in
B(K,\,H)}\sigma_{SF-}(M_{C}))\cup \{\lambda\in {\mathbb{C}}: A-\lambda\, \makebox{ is compact}\}.$$

(ii). Suppose that $\sigma_{*}(\cdot)=\sigma_{r}(\cdot)$ and 
$\lambda \not\in (\bigcap_{C\in B(K,\,H)}\sigma_{r}(M_{C}))\cup
\{\lambda\in {\mathbb{C}}: A-\lambda\, \makebox{ is compact}\}.$
 Then $A-\lambda$ is not compact and there exists $C\in B(K,\,H)$ such
 that $M_{C}-\lambda\in G_r(H\oplus K),$  and hence  $B-\lambda\in G_r(K).$

Case 1. $ \alpha(B-\lambda) =\infty$. It follows from Lemma 6 
that there exists $C\in \Phi(K, H)$ such that $M_C-\lambda$ is a lower 
semi-Weyl operator and $\beta(M_C-\lambda)=\beta(B-\lambda)$.
Note that $B-\lambda$ is surjective, then $M_C-\lambda$ is also
surjective. This implies that $\lambda \not\in \bigcap_{C\in
\Phi(K,\,H)}\sigma_{r}(M_{C})$.  It is clear that $$\bigcap_{C\in \Phi(K,\,H)}
\sigma_{r}(M_{C})\subseteq \bigcap_{C\in B(K,\,H)}\sigma_{r}(M_{C})\cup
\{\lambda\in{\mathbb{C}}: A-\lambda\, \makebox{ is compact}\}.$$

Case 2. $\alpha(B-\lambda)<\infty$. This means that $B-\lambda\in\Phi(K),$
so it is easy to prove that  $A-\lambda\in\Phi_-(H).$  Moreover, 
it follows from [10, Corollary 2] that $\alpha (B-\lambda)\geq \beta (A-\lambda).$
Next we claim that there exists some $C\in \Phi(K,\,H)$ such that 
$\lambda \not\in \bigcap_{C\in \Phi(K,\,H)}\sigma_{r}(M_{C})$. 
For this, let ${N(B-\lambda)}^\perp= K_1 \oplus K_2$ with $\dim K_2
=\dim \beta(A-\lambda)$. Define an operator $Q:
K\rightarrow H$ by $$Q=\left(\begin{array}{cc}C_1&0\\0&C_2\\\end{array}
 \right):({N(B-\lambda)}\oplus K_1)\oplus K_2\longrightarrow R(A-\lambda)\oplus R(A-\lambda)^\perp,$$
where $C_1\in B({N(B-\lambda)}\oplus K_1, R(A-\lambda))$ and $C_2\in B(K_2, R(A-\lambda)^\perp)$ are invertible operators. 
Obviously, operator $Q\in G(K,\,H)$ and $M_C-\lambda$ is  surjective. 
Thus $$\bigcap_{C\in G(K,\,H)}\sigma_{r}(M_{C})
\subseteq \bigcap_{C\in B(K,\,H)}\sigma_{r}(M_{C})\cup \{\lambda\in
{\mathbb{C}}: A-\lambda\, \makebox{ is compact}\}.$$

Together Case 1 with Case 2, we have 

$$\bigcap_{C\in\Phi(K,\,H)}\sigma_{r}(M_{C})=\bigcap_{C\in
B(K,\,H)}\sigma_{r}(M_{C})\cup \{\lambda\in {\mathbb{C}}:A-\lambda\, \makebox{ is compact}\}.$$

Similarly, when $\sigma_{*}=\sigma_{sw},$ we can prove the conclusion is also true.

\vspace{3mm}

By the proof methods of Theorem 4, we can prove the following result:

\vspace{3mm}
\noindent{\bf Theorem 5.} For a given pair $(A, B)\in B(H)\times B(K)$, we have 
$$\bigcap_{C\in\Phi(K,\,H)}\sigma_{*}(M_{C})=(\bigcap_{C\in B(K,\,H)}\sigma_{*}(M_{C}))\cup \{\lambda\in
 {\mathbb{C}}:B-\lambda\, \makebox{is   compact}\},$$ where
$\sigma_{*} \in\{\sigma_{l}, \sigma_{aw}\}.$


\begin{thebibliography}{S2}

\bibitem{1}  X. H. Cao. Browder spectra for upper triangular operator matrices, J. Math. Anal. Appl., 342(2008), 477-484.
\bibitem{2}  X. H. Cao, M. Z. Guo,  B. Meng. Semi-Fredholm spectrum and Weyl's theory for operator matrices, Acta Math. Sinica, 22(2006), 169-178.
\bibitem{3}  X. H. Cao, B. Meng. Essential approximate point spectra and Weyl's theorem for upper triangular operator matrices, J. Math. Anal. Appl., 304(2005), 759-771.
\bibitem{4}  X. L. Chen, S. F. Zhang, H. J. Zhong. On the filling in holes problem for operator matrices, Linear Algebra Appl., 430(2009),558-563
\bibitem{5}  D. S. Djordjevi\'{c}. Perturbations  of spectra of operator matrices, J. Operator Theory, 48(2002), 467-486.
\bibitem{6}  S. V. Djordjevi\'{c}, Y. M. Han. spectral continuity for operator matrices, Glasg. Math. J., 43(2001), 487-490.
\bibitem{7}  S. V. Djordjevi\'{c},  H. Zguitti. Essential point spectra of operator matrices though local spectral theory, J.  Math. Anal. Appl., 338(2008), 285-291.
\bibitem{8}  H. K. Du, J. Pan. Perturbation of  spectrums  of $2\times2$ operator matrices,  Proc. Amer. Math. Soc., 121(1994), 761-766.
\bibitem{9}  J. K. Han, H. Y. Lee, W. Y. Lee. Invertible completions of $2\times2$ upper triangular operator matrices. Proc. Amer. Math. Soc., 128(1999), 119-123.
\bibitem{10} I. S. Hwang, W. Y.  Lee. The boundedness below of $2\times2$ upper triangular operator matrices, Integr. Equ. Oper. Theory, 39(2001), 267-276.
\bibitem{11} W. Y. Lee. Weyl spectra of operator matrices, Proc. Amer. Math. Soc., 129(2000), 131-138.
\bibitem{12} Y. Li, H. K. Du. The intersection of left and right essential spectra of 2 $\times$ 2 operator matrices, Bull. Lond. Math. Soc. , 36(2004), 811-819.
\bibitem{13} Y. Li, H. K. Du. The intersection of essential approximate point spectra of operator matrices, J. Math. Anal.  Appl., 323(2006), 1171-1183.
\bibitem{14} Y. Li, X. H. Sun, H. K. Du. The intersection of left(right) spectra  of 2 $\times$ 2 upper triangular operator matrices, Linear Algebra Appl., 418(2006), 112-121.
\bibitem{15} Y. Li, X. H. Sun, H. K. Du. A note on the left essential approximate point spectra  of operator matrices, Acta Math. Sinica,  23(2007), 2235-2240.
\bibitem{16} E. H. Zerouali,  H. Zguitti. Perturbation of spectra of operator matrices and local spectral theory, J. Math. Anal. Appl., 324(2006), 992-1005.
\bibitem{17} S. F. Zhang, H. J. Zhong. A note of Browder spectrum of operator matrices, J. Math. Anal. Appl., 344(2008), 927-931.
\bibitem{18} S. F. Zhang, H. J. Zhong, Q. F. Jiang.  Drazin  spectrum of operator matrices on the Banach space, Linear Algebra Appl., 429(2008), 2067-2075.
\bibitem{19} Y. N. Zhang, H. J. Zhong, L. Q. Lin.  Browder spectra and essential spectra of operator matrices, Acta Math. Sinica, 24(2008), 947-954.

\end{thebibliography}
\end{document}